\documentclass{amsart}
\usepackage{amsmath, amsthm, amscd, amsfonts, amssymb, graphicx, color}

\usepackage{longtable}
\usepackage{hyperref} 
\usepackage{cleveref} 
\usepackage{listings}
\usepackage{algorithm}
\usepackage{algpseudocode}
\makeatletter
\renewcommand{\ALG@name}{Sub-Algorithm}
\makeatother
\makeatletter
\newenvironment{breakablealgorithm}
{
		\begin{center}
			\refstepcounter{algorithm}
			\hrule height.8pt depth0pt \kern2pt
			\renewcommand{\caption}[2][\relax]{
				{\raggedright\textbf{\ALG@name~\thealgorithm} ##2\par}%
				\ifx\relax##1\relax 
				\addcontentsline{loa}{algorithm}{\protect\numberline{\thealgorithm}##2}%
				\else 
				\addcontentsline{loa}{algorithm}{\protect\numberline{\thealgorithm}##1}%
				\fi
				\kern2pt\hrule\kern2pt
			}
		}{
		\kern2pt\hrule\relax
	\end{center}
}
\makeatother

\newtheorem{theorem}{Theorem}[section]
\newtheorem{lem}[theorem]{Lemma}

\theoremstyle{definition}

\newtheorem{conj}[theorem]{Conjecture}
\newtheorem{prop}[theorem]{Proposition}

\theoremstyle{remark}

\numberwithin{equation}{section}

\newcommand{\gyr}[2]{{\mathrm{gyr}[{#1}]}{#2}}

\newcommand{\Aut}[1]{\mathrm{Aut}\,{(#1)}}

\begin{document}

\title[Construction of All Gyrogroups of Orders at most 31]{Construction of All Gyrogroups of Orders at most 31}

\author[A. R. Ashrafi]{ Ali Reza Ashrafi }
\address{Department of Pure Mathematics, Faculty of Mathematical Sciences, University of Kashan, Kashan 87317-53153, I. R. Iran (https://orcid.org/0000-0002-2858-0663)}
\email{ashrafi@kashanu.ac.ir}
\thanks{}

\author[K. Mavaddat Nezhaad]{Kurosh Mavaddat Nezhaad}
\address{Department of Pure Mathematics, Faculty of Mathematical Sciences, University of Kashan, Kashan 87317-53153, I. R. Iran (https://orcid.org/0000-0001-8260-2619)}
\email{kuroshmavaddat@gmail.com}
\thanks{}

\author[M. A. Salahshour]{Mohammad Ali Salahshour}
\address{Department of Mathematics, Savadkooh Branch, Islamic Azad University, Savadkooh, I. R. Iran (https://orcid.org/0000-0002-0816-4232)}
\email{salahshour@iausk.ac.ir}
\thanks{}


\date{}

\dedicatory{This work is dedicated to professor Abraham Ungar for his pioneering role in gyrostructures}

\begin{abstract}
The gyrogroup is the closest algebraic structure to the group ever discovered. It has a binary operation $\star$ containing an identity element such that each element has an inverse. Furthermore, for each pair $(a,b)$ of elements of this structure there exists an automorphism $\gyr{a,b}{}$ with this property that left associativity and left loop property are satisfied. Since each gyrogroup is a left Bol loop, some results of Burn imply that all gyrogroups of orders $p, 2p$ and $p^2$ are groups. The aim of this paper is to classify gyrogroups of orders 8, 12, 15, 18, 20, 21, and 28.
\end{abstract}

\maketitle

\section{Introduction}
At the dawn of 20th century, Albert Einstein changed the face of modern physics and astronomy, by publishing four papers in 1905, which also known as annus mirabilis papers \cite{5}.  In one of these papers, he was  introduced a new addition $\oplus$ on $\mathbb{R}^3$ such that $(\mathbb{R}^3,\oplus)$ has identity and each element of $\mathbb{R}^3$ has an inverse, but the associativity is not satisfied. Instead of associativity it satisfies another condition which is now called  gyroassociaive law. The theory of gyrostructures is a result of some innovating idea developed by Abraham Ungar during the ninth decade of the 20th century \cite{8,9}.

In the first half of twentieth century by the impressions of geometry and algebra, Albert \cite{1,2}, Baer \cite{3} with more interest in relation with geometry and Bruck \cite{4} with more interest in relation with algebra contributed to the theory of loops and quasigroups which later appeared as an interpretation of Ungar's gyrogroups \cite{8}.

In 1989 Abraham Ungar  continued his work  to discover a pattern behind the seemingly lawless Einstein's addition of velocities \cite{10}. He founded gyrogroup theory, non-associative groups that share analogies with groups, and also provides a context for analytic hyperbolic geometry. In the past years, Ungar and the others proved many important results in gyrogroup theory and hyperbolic geometry, which have applications in foundations of physics.

We are now ready to present the Ungar's definition of gyrogroup directly. According to the Ungar's famous book \cite{11,12},  a pair $(G,\oplus)$ consists of a nonempty set $G$ together with a binary operation $\oplus$ on $G$ is called a \textit{gyrogroup} if its binary operation satisfies the following axioms:
\begin{enumerate}
\item [(G1)] there exists an element $0 \in G$ such that for all $x \in G$, $0 \oplus x = x$;

\item [(G2)] for each $a \in  G$, there exists $b \in G$ such that $b \oplus a$ $= 0$;

\item [(G3)] there exists a function $\mbox{gyr}: G \times G \longrightarrow \Aut{G}$ such that for every $a, b, c \in G$, $a \oplus (b \oplus c) =  (a \oplus b) \oplus \gyr{a,b}{c}$;

\item [(G4)] for each $a, b \in G$, $\gyr{a,b}{} = \gyr{a \oplus b, b}{}$.
\end{enumerate}

The map $\gyr{a, b}{}$ is called the \textit{gyroautomorphism 	generated by $a$ and $b$}. A gyrogroup $G$ can have other property like gyrocommutativity that can be formulated as $a\oplus b = \gyr{a, b}({b\oplus a})$, for all $a, b\in G$. A non-degenerate gyrogroup is a gyrogroup with a nontrivial gyroautomorphism. It is easy to see that a gyrogroup $G$ is a group if and only if  $G$ does not have non-trivial gyroautomorphism. This shows that a non-degenerate gyrogroup cannot be a group.

For the sake of completeness, we mention here a result which is crucial throughout this paper. The interested readers can consult \cite[pp. 19-20]{11} for its proof.

\begin{theorem}\label{th1}
A magma $(G,\oplus)$ forms a gyrogroup if and only if it satisfies the following conditions:
\begin{enumerate}
\item  there exists $0\in G$ such that for all $a\in G$, we have $0\oplus a = a = a\oplus 0$;

\item for each $a\in G$ there exists $b\in G$ such that $b\oplus a = 0 = a\oplus b$;

\item if we define $\gyr{a,b}{c} = \ominus(a\oplus	b)\oplus(a\oplus (b\oplus c))$, then $\gyr{a,b} \in\Aut{G,\oplus}$;

\item $a\oplus (b\oplus c) = (a\oplus b)\oplus\gyr{a, b}{c}$;

\item $(a\oplus b)\oplus c = a\oplus (b\oplus\gyr{b, a}{c})$;

\item $\gyr{a, b}{} = \gyr{a\oplus b, b}{}$;

\item $\gyr{a, b}{} = \gyr{a, b\oplus a}{}$.
\end{enumerate}
\end{theorem}\

By Theorem \ref{th1}, any gyroautomorphism is completely determined by its generators via the \textit{gyrator identity}.

Ungar discovered that the M$\ddot{\rm o}$bius gyrogroup is a gyrocommutative gyrogroup and the set of all gyrators is not a subgroup of the whole automorphism group of the M$\ddot{\rm o}$bius gyrogroup. In all known examples of finite gyrogroups, the set of all gyrators of a finite gyrogroup forms a subgroup of the full automorphism group.

The main results of this paper are as follows:

\begin{theorem} \label{tt1}
There exists a gyrogroup of order 16 in which the set of all gyrators is not a group.
\end{theorem}

\begin{theorem} \label{tt2}
Up to isomorphism, there are six non-degenerate gyrogroups of order 8; two non-degenerate gyrogroups of  orders 12, 20 and 28; one non-degenerate gyrogroup of both orders 15 and 21. There is no non-degenerate gyrogroups of orders a prime number, a prime square, two times of a prime number and 18. See Table \ref{ta2} for details.
\end{theorem}

\section{Classification of Non-Degenerate Gyrogroups of Orders $\leq 31$}
The aim of this section is to present the classification of gyrogroups of orders $\leq 31$ except from the orders $24, 27$ and $30$. Our classifications are based on Burn's results \cite{42,44,46}  and calculations given by  Moorhouse \cite{55}.

\subsection{An Algorithm for constructing gyrogroups from Bol loops.}

Sabinin et al. \cite{6} proved that all gyrogroups are left Bol loops and Burn \cite{42} proved that all Bol loops of orders $p$, $2p$ and $p^2$, $p$ is prime, are groups. So, it is enough to consider all left Bol loops of a given order and then check the properties of a gyrogroup. By these results we prepare an algorithm for constructing gyrogroups of a given order $n$ based on the structure of all left Bol loops of order $n$. To see our algorithm, we choose a left Bol loop $L$ of order $n$. For each $a, b \in L$, we define the mapping $\gyr{a,b}{}: L \longrightarrow L$ by $\gyr{a,b}(x) = \ominus(a\oplus	 b)\oplus(a\oplus (b\oplus x))$. If all mappings $\gyr{a,b}{}$ are automorphisms of $L$ and these automorphisms satisfy the condition $G3$ then $L$ will be a gyrogroup. Note that by \cite[p. 13]{6} the left loop property (G4) is equivalent to the left Bol identity and so we don't need to check left loop property.

To construct all gyrogroups of a given order, it is enough to choose a left Bol loop $K$ with binary operation $\oplus$. By Theorem \ref{th1}(3),  all gyrogroups satisfy the equation $\gyr{a,b}{c} = \ominus(a\oplus	 b)\oplus(a\oplus (b\oplus c))$, where $a, b$ and $c$ are arbitrary elements of $K$. This implies that, if $K$ is a gyrogroup then this equality should be satisfied. The following simple Gap code \cite{75} define  a function $\gyr{ \ }{ }$ based on three elements $a, b$ and $x$ from $B$.

\vskip 3mm

\begin{breakablealgorithm}
\caption{Computing gyr[a,b](x)}\label{alg:cap}
\begin{algorithmic}
	\State\textbf{Input:} Elements $a$ and $b$ of a loop $B$;
	\State\textbf{Output:} The function $\gyr{a,b}{}$ by definition given Theorem \ref{th1};\\
	\Function{gyr}{}\\
	\hspace{.5cm}\Return $((a*b)^{-1})*(a*(b*x))$;\Comment{$g(a,b)(x)=-(a \oplus b) \oplus (a \oplus (b \oplus x))$}
	\EndFunction
\end{algorithmic}
\end{breakablealgorithm}

We now assume that $A$ is the set of all gyrators $\gyr{a,b}{}$, where $a$ and $b$ are elements of the left Bol loop $K$.  Again, if $K$ is a gyrogroup then Theorem \cite[Theorem 1.13]{12} shows that for each element $b \in B$,  $\gyr{a,a}{}$ $=$ $\gyr{a,-a}{}$ $=$ $\gyr{0,a}{}$ $=$ $\gyr{a,0}{}$ $=$ $I$. Therefore, if one of these equalities are not satisfied, then $K$ will not a gyrogroup.

\begin{breakablealgorithm}
	\caption{Computing the Non-Identity Gyroautomorphisms}\label{alg:cap}
	\begin{algorithmic}
		\State\textbf{Input:} Two elements $a$ and $b$ of $K$.
		\State\textbf{Output:} Add $\gyr{a,b}{}$ to $A$, when $\gyr{a,b}{}$ is a non-identity automorphism.
		\Function{gyroauto}{}\\
		\hspace{.5cm}local m,a,b,x,y,z,A,B,C,D,DD,EK,EE,FF,KK
		\State A:=[ ]
		\State KK:= Difference(K,[Identity(K)])\Comment{text}
		\For{$a \in KK$} \do\\
		\For{$b \in KK$} \do\\
		\If{$ a\neq b$ and $b\neq a^{-1}$} \\\Comment{text}
		\State B:=[ ]\\
		\State C:=[ ]\\
		\For{$x \in K$} \do\\
		\State Add(C,x)\\
		\State Add(B,gyr(a,b,x))\\
		\EndFor
		\If{$B\neq C$}  \Comment{g(a,b) is non-identity permutation}\\
		\State Append(A,[[a,b],B])\\
		\EndIf\\
		\EndIf\\
		\EndFor\\
		\EndFor\\
		\State m:=Size(A)\\
		\State D:=Set(A{[2,4..m]})\\
		\State EE:=[ ]\\
		\For{$y \in D$} \do\\
		\State FF:=[ ]\\
		\For{$z \in y$} \do\\
		\State Add(FF,Position(K,z))\\
		\EndFor\\
		\State Add(EE,FF)\\
		\EndFor\\
		\State DD:=List(EE,x$\rightarrow$PermList(x))\\
		\If{IsSubset(AutomorphismGroup(K),DD)} \\
		\State \Return A\\
		\Else
		\State \Return false\\
		\EndIf\\
		\EndFunction\\
	\end{algorithmic}
\end{breakablealgorithm}

By a result of Sabinin \cite{6}, all mappings $\gyr{a,b}{}$ are bijective and so it is enough to check that whether or not $\gyr{a,b}{}$ is an automorphism.

\begin{breakablealgorithm}
	\caption{Computing Gyroautomorphism Table of a Gyrogroup}
	\begin{algorithmic}
		\Function{matgyroauto}{}
		\State local a,b,i,j,s,ss,t,n,m,x,y,z,M,B,D,DD,DDD,EE,FF
		\State M:=[ ]
		\State B:=[``A",``B",``C",``D",``E",``F",``G",``H",``K",``L",``M",``N",``P",``Q",``R",

                 \hspace{0.5cm}  ``S",``T",``U",``V",``W",``X",``Y",``Z",``J",``O"]
		\State m:=Size(A)
		\State D:=A{[1,3..m-1]}
		\State DD:=Set(A{[2,4..m]})
		\State EE:=[ ]
		\For{$y \in DD$} \do\\
		\State FF:=[ ]
		\For{$z \in y$} \do\\
		\State Add(FF,Position(K,z))\\
		\EndFor\\
		\State Add(EE,FF)
		\EndFor\\
		\State DDD:=List(EE,x$\longrightarrow$PermList(x))\\
		\For{$a \in K$} \do\\
		\State i:=Position(K,a)
		\State M[i]:=[ ]
		\For{$b \in K$} \do\\
		\State j:=Position(K,b)
		\If{ $[a,b] \in D$}
		\State t:=Position(A,[a,b])
		\State s:=A[t+1]
		\State ss:=Position(DD,s)
		\State M[i][j]:=B[ss]
		\Else\\
		\State M[i][j]:=``I"
		\EndIf\\
		\EndFor\\
		\EndFor\\
		\State Print(``Non-identity automorphisms are as follows:",DDD)
		\State \Return M
		\EndFunction
	\end{algorithmic}
\end{breakablealgorithm}

Sub-Algorithm 4 determined that whether or not a left Bol loop $K$ is a gyrogroup.

\begin{breakablealgorithm}
	\caption{A function which shows that $M$ is the Cayley table of a gyrogroup}
	\begin{algorithmic}
		\Function{isgyro}{}
		\State local k,a,b,c,n,A,B,N,K,KK
		\If{IsLoopTable(M)=false} \State \Return false \EndIf\Comment{Check Loop Table}
		\State K:=LoopByCayleyTable(M)
		\State KK:= Difference(K,[Identity(K)])
		\For{$a \in KK$} \do\\
		\For{$b \in KK$} \do\\
		\If{$b\neq a^{-1}$}
		\For{$c \in KK$} \do\\
		\If{$a*(b*c)\neq (a*b)*\gyr{a,b}{c}$}
		\State \Return false
		\EndIf \Comment{Check Gyroassociativity}\\
		\EndFor\\
		\EndIf
		\EndFor\\
		\EndFor\\
		\State N:=gyroauto(K)\Comment{Check the automorphism conditions for gyrators of M}
		\If{N=false}
		\State \Return false
		\ElsIf N=[ ]
		\State Print(``This is a group")
		\Else
		\State Print(``Cayley table and gyroautomorphism table")
		\State Print(matgyroauto(K,N))
		\EndIf
		\EndFunction
	\end{algorithmic}
\end{breakablealgorithm}
\par

We end this subsection by introducing an example of a gyrogroup of order 16 in which the set of all gyroautomorphisms is not a group. To see this, we assume that $G_{16}$ is a left Bol loop of order 16 with the Cayley table \ref{ta1}:

\begin{table}\caption{The Cayley Table of $G_{16}$}\label{ta1}
\begin{center}
\begin{tabular}{c|ccccccccccccccccc}
        $\oplus$ & 0 & 1 & 2 & 3 & 4 & 5 & 6 & 7 & 8 & 9 & 10 & 11 & 12 & 13 & 14 & 15 \\ \hline
		0 & 0 & 1 & 2 & 3 & 4 & 5 & 6 & 7 & 8 & 9 & 10 & 11 & 12 & 13 & 14 & 15 \\
		1 & 1 & 0 & 4 & 6 & 2 & 7 & 3 & 5 & 10 & 15 & 8 & 12 & 11 & 14 & 13 & 9 \\
		2 & 2 & 3 & 5 & 7 & 1 & 6 & 0 & 4 & 12 & 11 & 14 & 10 & 15 & 8 & 9 & 13 \\
		3 & 3 & 2 & 7 & 5 & 0 & 4 & 1 & 6 & 11 & 12 & 13 & 15 & 10 & 9 & 8 & 14 \\
		4 & 4 & 6 & 1 & 0 & 5 & 3 & 7 & 2 & 14 & 13 & 12 & 8 & 9 & 10 & 15 & 11 \\
		5 & 5 & 7 & 6 & 4 & 3 & 0 & 2 & 1 & 15 & 10 & 9 & 14 & 13 & 12 & 11 & 8 \\
		6 & 6 & 4 & 0 & 1 & 7 & 2 & 5 & 3 & 13 & 14 & 11 & 9 & 8 & 15 & 10 & 12 \\
		7 & 7 & 5 & 3 & 2 & 6 & 1 & 4 & 0 & 9 & 8 & 15 & 13 & 14 & 11 & 12 & 10 \\
		8 & 8 & 9 & 11 & 12 & 13 & 15 & 14 & 10 & 5 & 7 & 1 & 6 & 4 & 3 & 2 & 0 \\
		9 & 9 & 8 & 12 & 11 & 14 & 10 & 13 & 15 & 7 & 5 & 0 & 4 & 6 & 2 & 3 & 1 \\
		10 & 10 & 15 & 13 & 14 & 11 & 9 & 12 & 8 & 1 & 0 & 5 & 3 & 2 & 6 & 4 & 7 \\
		11 & 11 & 12 & 8 & 9 & 10 & 14 & 15 & 13 & 2 & 3 & 4 & 0 & 1 & 7 & 5 & 6 \\
		12 & 12 & 11 & 9 & 8 & 15 & 13 & 10 & 14 & 3 & 2 & 6 & 1 & 0 & 5 & 7 & 4 \\
		13 & 13 & 14 & 15 & 10 & 9 & 12 & 8 & 11 & 6 & 4 & 3 & 7 & 5 & 0 & 1 & 2 \\
		14 & 14 & 13 & 10 & 15 & 8 & 11 & 9 & 12 & 4 & 6 & 2 & 5 & 7 & 1 & 0 & 3 \\
		15 & 15 & 10 & 14 & 13 & 12 & 8 & 11 & 9 & 0 & 1 & 7 & 2 & 3 & 4 & 6 & 5 \\
\end{tabular}
\end{center}
\end{table}

Our Gap code shows that $G_{16}$ is a gyrogroup  and  all non-identity gyroautomorphisms of $G_{16}$ are as follows:
\begin{center}
\begin{tabular}{ll}
$(3,4)(5,7)(9,10)(11,16)$ & $(3,5)(4,7)(9,11)(10,16)$\\
$(3,5)(4,7)(9,16)(10,11)(12,13)(14,15)$ & $(3,7)(4,5)(9,11)(10,16)(12,13)(14,15)$\\
$(3,7)(4,5)(9,16)(10,11)$ &
\end{tabular}
\end{center}

It is an elementary fact that there is no elementary abelian group of order six and so we don't have a subgroup of $Aut(G_{16})$. This proves Theorem \ref{tt1}.

\subsection{Gyrogroups of Orders 8, 12, 16, 18, 20, 21, 28}

The aim of this subsection is to given a classification of gyrogroups of orders 8, 12, 16, 18, 20, 21, 28. Some gyrogroups of orders 27 are also given.

Mahdavi et al. \cite{53} asked about the number of gyrogroups of order 8. The following lemma respond to this question.

\begin{lem}
There are exactly six gyrogroups of order 8.
\end{lem}

\begin{proof}
Burn \cite[Theorem 6]{42} proved that there are exactly six non-associative Bol loops of order 8. So, it is enough to check these Bol loops one by one. By our Gap code given in Subsection 2.1, it can be shown that all Bol loops of order 8 are gyrogroups and so there are exactly six gyrogroups of order 8.
\end{proof}

\vskip 3mm

\begin{lem}
There are exactly two gyrogroups of each order 12, 20 and 28.
\end{lem}

\begin{proof}
By Burn \cite[Theorems 1 and 2]{44}, there are exactly three non-associative Bol loops of order $4p$, $p$ is prime, such that only one of them is a Moufang loop. Our calculations by the Gap code  in Subsection 2.1 shows that the Moufang loop is not a gyrogroup, but both Bol loops which are not of Moufang type are gyrogroups. Therefore, there are exactly two gyrogroups of each order 12, 20 and 28.
\end{proof}

\vskip 3mm

\begin{lem}
There is a unique gyrogroup of order 15 which is gyrocommutative.
\end{lem}

\begin{proof}
Niederreiter and  Robinson \cite{57} investigated the structure of Bol loops of order $pq$, where $p$ and $q$ are different prime and $q$ divides $p^2 - 1$. As a consequence, there are exactly two Bol loops of order 15. Our Gap code shows that exactly one of these Bol loops of order 15 is a gyrogroup which is the gyrogroup reported by Suksumran  in \cite[p. 432]{73}.
\end{proof}

\vskip 3mm

\begin{lem}
There is no gyrogroup of order 18.
\end{lem}

\begin{proof}
Burn \cite[Theorem 6]{46, 48}, proved that there are exactly two non-associative left Bol loop of order $2p^2$ which is not a Moufang loop. Our Gap code shows that for the case that $p = 3$, these Bol loops are not a gyrogroup. So, there is no gyrogroup of order 18.
\end{proof}

\vskip 3mm

\begin{lem}
There is a unique gyrogroup of order 21 which is gyrocommutative.
\end{lem}

\begin{proof}
Kinyon et al. \cite[Theorem 1.1]{51}, proved that if $p$ and $q$ are primes such that $q$ divides $p^2-1$, then there exists a unique non-associative left Bruck loop of order $pq$, up to isomorphism, and there are precisely $\frac{(p -q+4)}{2}$ left Bol loops of order $pq$. This proves that there are exactly four Bol loops of order 21 that two of them are groups. Again we apply our Gap code to prove that precisely one of these non-associative Bol loop is a gyrogroup.
\end{proof}

\vskip 3mm

\noindent{\bf Proof of Theorem \ref{tt2}}. By Lemmas 2.1 - 2.5, there are exactly six non-degenerate gyrogroups of order 8; two non-degenerate gyrogroups of each  order 12, 20 and 28; one non-degenerate gyrogroup of both orders 15 and 21. There is no non-degenerate gyrogroups of orders a prime number, a prime square, two times of a prime number and 18.

\vskip 3mm

There is no classification of Bol loops of orders 24, 27 and 30 and since our algorithm is based on the classification of Bol loops, we don't have all gyrogroups of these orders, but we can construct 8 gyrogroups of order 28 and one gyrogroup of order 30. So, we don't have information about gyrogroups of order 24. Finally, it is possible to construct 1995 gyrogroups of order 16.  We record our results in Table  \ref{ta2}.  In this table, $\alpha(n)$ and $\beta(n)$ denote the number of all non-degenerate gyrogroups and non-degenerate gyrocommutative gyrogroups of order $n$, respectively.

\begin{table}
\caption{The Number of All Non-Degenerate Gyrogroups and Non-Degenerate Gyrocommutative Gyrogroups}\label{ta2}
\begin{tabular}{c|c|c}
n & $\alpha(n)$  & $\beta(n)$ \\  \hline
		8 & 6 & 3 \\
		12 & 2 & 0 \\
		15 & 1 & 1 \\
		16 & 1995 & 179 \\
		18 & 0 & 0 \\
		20 & 2 & 0 \\
		21 & 1 & 1 \\
		24 & - & - \\
		27 & $\geq$8 & $\geq$4 \\
		28 & 2 & 0 \\
		30 & $\geq$1 & $\geq$1\\
\end{tabular}
\end{table}

\section{Concluding Remarks}
In this paper, a classification of gyrogroups of orders less than 32 except from 24, 27 and 30 are given. Our argument are based on the Burn's classification of Bol loops \cite{42,44,46} and calculations of Moorhouse \cite{55}. Eight gyrogroups of orders 27 and one gyrogroup of order 30 are also constructed. All tables of these gyrogroups can be downloaded from https://faculty.kashanu.ac.ir/ashrafi/en.

Suppose $G$ is a gyrogroup and $a, b \in G$. We define the commutator of $a$ and $b$ as $[a,b]$ $=$ $\ominus(a \oplus b) \oplus \gyr{a,b}{(b \oplus a)}$. The subgyrogroup generated by all commutators of $G$, $G^\prime$, is called the derived subgyrogroup of $G$. Suksumran \cite{7} proved that $G^\prime$ is a subgroup of $G$ and conjectured that it is normal subgroup. We check this conjecture on all gyrogroups of orders  less than 32 except from orders 24, 27 and 30 and we have the following result:

\begin{prop}
The commutator subgroup of all gyrogroups of orders  less than 32 except from orders 24, 27 and 30 are normal in the whole gyrogroup.
\end{prop}

\vskip 3mm

Our calculations with the aid of Gap suggests the following conjecture on finite gyrogroups:

\vskip 3mm

\begin{conj}
If a gyrogroup is commutative then it is group.
\end{conj}

\vskip 3mm

\noindent\textbf{Data Availability} The data used to support  this study can be downloaded from https://faculty.kashanu.ac.ir/ashrafi/en.

\vskip 3mm

\noindent\textbf{Conflicts of Interest} The authors declare that they have no conflicts of interest.

\vskip 3mm

\noindent\textbf{Authorship Statement:} 
\textbf{Ali Reza Ashrafi:} Conceptualization, Methodology, Writing-Original draft preparation.
\textbf{Kurosh Mavaddat Nezhaad:} Software, Data curation, Investigation, Editing.
\textbf{Mohammad Ali Salahshour:} Software, Investigation, Validation.

\vskip 3mm

\bibliographystyle{amsplain}

\end{document}